\documentclass[a4paper,10pt]{article}
\usepackage[dvips]{epsfig}
\usepackage{amsthm,amsmath,amssymb}
\usepackage{latexsym}
\usepackage{graphics}
\usepackage{url}
\usepackage{placeins}
\usepackage{hyperref}
\usepackage{color}
\usepackage{bm}
\usepackage{comment}
\usepackage{statex}

\definecolor{blue}{rgb}{0,0,0.9}
\definecolor{red}{rgb}{0.9,0,0}
\definecolor{green}{rgb}{0,0.50,0.10}
\definecolor{violet}{rgb}{0.5804,0.0000,0.8275}  

\newcommand{\blue}[1]{\begin{color}{blue}#1\end{color}}

\def\@themcountersep{}

\newtheorem{THEO}{Theorem}[section]
\newtheorem{ALGo}[THEO]{Algorithm}
\newtheorem{CONJ}[THEO]{Conjecture}
\newtheorem{COND}[THEO]{Condition}
\newtheorem{CORO}[THEO]{Corollary}
\newtheorem{DEFI}[THEO]{Definition}
\newtheorem{EXAMP}[THEO]{Example}
\newtheorem{FACT}[THEO]{Fact}
\newtheorem{HYPO}[THEO]{Hypothesis}
\newtheorem{LEMM}[THEO]{Lemma}
\newtheorem{PROB}[THEO]{Problem}
\newtheorem{PROP}[THEO]{Proposition}
\newtheorem{REMA}[THEO]{Remark}
\newcommand{\theo}{\begin{THEO}}
\newcommand{\algo}{\begin{ALGo} \rm}
\newcommand{\cond}{\begin{COND}}
\newcommand{\conj}{\begin{CONJ}}
\newcommand{\coro}{\begin{CORO}}
\newcommand{\defi}{\begin{DEFI} \rm}
\newcommand{\examp}{\begin{EXAMP} \rm}
\newcommand{\fact}{\begin{FACT}}
\newcommand{\hypo}{\begin{HYPO} \rm}
\newcommand{\lemm}{\begin{LEMM}}
\newcommand{\prob}{\begin{PROB} \rm}
\newcommand{\prop}{\begin{PROP}}
\newcommand{\rema}{\begin{REMA} \rm}
\newcommand{\etheo}{\end{THEO}}
\newcommand{\ealgo}{\end{ALGo}}
\newcommand{\econd}{\end{COND}}
\newcommand{\econj}{\end{CONJ}}
\newcommand{\ecoro}{\end{CORO}}
\newcommand{\edefi}{\end{DEFI}}
\newcommand{\eexamp}{\end{EXAMP}}
\newcommand{\efact}{\end{FACT}}
\newcommand{\ehypo}{\end{HYPO}}
\newcommand{\elemm}{\end{LEMM}}
\newcommand{\eprob}{\end{PROB}}
\newcommand{\eprop}{\end{PROP}}
\newcommand{\erema}{\end{REMA}}

\def\0{\mbox{\bf 0}}
\def\1{\mbox{\bf 1}}
\def\2{\mbox{\bf 2}}
\def\3{\mbox{\bf 3}}
\def\4{\mbox{\bf 4}}
\def\5{\mbox{\bf 5}}
\def\6{\mbox{\bf 6}}
\def\7{\mbox{\bf 7}}
\def\8{\mbox{\bf 8}}
\def\9{\mbox{\bf 9}}
\def\a{\mbox{\boldmath $a$}}

\def\cc{\mbox{\boldmath $c$}}

\def\e{\mbox{\boldmath $e$}}

\def\s{\mbox{\boldmath $s$}}

\def\u{\mbox{\boldmath $u$}}

\def\x{\mbox{\boldmath $x$}}

\def\z{\mbox{\boldmath $z$}}
\def\A{\mbox{\boldmath $A$}}
\def\B{\mbox{\boldmath $B$}}
\def\C{\mbox{\boldmath $C$}}

\def\H{\mbox{\boldmath $H$}}

\def\O{\mbox{\boldmath $O$}}
\def\P{\mbox{\boldmath $P$}}
\def\Q{\mbox{\boldmath $Q$}}

\def\U{\mbox{\boldmath $U$}}

\def\X{\mbox{\boldmath $X$}}
\def\Y{\mbox{\boldmath $Y$}}

\def\EC{\mbox{$\cal E$}}

\def\NC{\mbox{$\cal N$}}

\def\inprod#1#2{\langle#1, \, #2\rangle}
\def\Inprod#1#2{\big\langle#1, \, #2\big\rangle}

\def\Real{\mbox{$\mathbb{R}$}}

\def\balpha{\mbox{\boldmath $\alpha$}}
\def\bbeta{\mbox{\boldmath $\beta$}}

\def\salpha{\mbox{\scriptsize \boldmath $\alpha$}}
\def\sbeta{\mbox{\scriptsize \boldmath $\beta$}}

\def\s0{\mbox{\scriptsize \boldmath $0$}}

\def\Real{\mathbb{R}}
\def\coneK{\mathbb{K}}

\def\SymMat{\mathbb{S}}

\def\lb{\mbox{lb}}
\def\ub{\mbox{ub}}
\def\BB{B\&B }
\def\vec{\mbox{\bf vec}}

\def\bQAP{\mbox{$\overline{\mbox{QAP}}$}}
\begin{document}

\title{ \Large 
Solving Challenging Large Scale QAPs  
} 

\author{
\normalsize 
Koichi Fujii\thanks{NTT DATA Mathematical Systems Inc., Tokyo, 160-0016 Japan}, \and \normalsize
Naoki Ito\thanks{Fast Retailing Co., Ltd., Uniqlo City Tokyo, Tokyo, 135-0063 Japan}, \and \normalsize
Sunyoung Kim\thanks{Department of Mathematics, Ewha W. University, 52 Ewhayeodae-gil, Sudaemoon-gu, Seoul 120-750, Korea 
}, \and \normalsize
Masakazu Kojima\thanks{Department of Industrial and Systems Engineering,
	Chuo University, Tokyo 192-0393, Japan 
},
  \and \normalsize
Yuji Shinano\thanks{Department of Applied Algorithmic Intelligence Methods (A$^2$IM), Zuse Institute Berlin, Takustr. 7, 14195 Berlin, Germany}, \and \normalsize
Kim-Chuan Toh\thanks{Department of Mathematics, and Institute of Operations Research and Analytics, National University of Singapore,
10 Lower Kent Ridge Road, Singapore 119076
} 
}

\date{\normalsize\today}
\maketitle 






\begin{abstract}

We report our  progress on the project for solving larger scale quadratic assignment problems (QAPs).
Our main approach to solve
 large scale NP-hard combinatorial optimization 
problems such as QAPs is a parallel 
branch-and-bound method efficiently implemented on a powerful computer system 
using the Ubiquity Generator (UG) framework that can utilize 
more than 100,000 cores.
 Lower bounding procedures incorporated in the branch-and-bound method play a crucial role in solving the problems.
 For a strong lower bounding procedure, 
we employ the Lagrangian doubly 
nonnegative (DNN) relaxation and the Newton-bracketing method developed 
by the authors' group. In this report, we describe 
some basic tools used in the project including the 
lower bounding procedure and branching rules, and present some preliminary numerical results.
 
Our next target problem is QAPs with dimension at least 50, as we 
have succeeded to solve tai30a and sko42 from QAPLIB for the first time.

\end{abstract}



\section{Introduction}

For a positive integer $n$, we let 
$N=\{1,\ldots,n\}$ represent a set of locations and also a set of facilities.
Given $n \times n$ symmetric matrices $\A =[a_{ik}]$, 
$\B=[b_{j\ell}]$ and an $n \times n$ matrix $\C  = [c_{ij}]$, 
the quadratic assignment problem (QAP) is 
stated as 
\begin{eqnarray}
\min_{\pi} \sum_{i \in N}^n\sum_{k\in N}^n a_{ik}b_{\pi(i),\pi(k)}
+\sum_{i\in N}^n c_{i,\pi(i)}, \label{eq:QAP0}
\end{eqnarray}
where $a_{ik}$ denotes the flow between facilities $i$ and $k$, $b_{j\ell}$ is 
the distance between locations $j$
and $\ell$, $c_{ij}$ the fixed cost of assigning facility $i$ to location $j$, and $(\pi(1),\ldots,\pi(n))$ a permutation of 
$1,\ldots,n$ such that 
$\pi(i) = j$ if facility $i$ is
assigned to location $j$. 

The QAP is known as NP-hard in theory, and
solving QAP instances of size larger than 35 in practice   still remains challenging. 
Various heuristic methods 
for the QAP such as tabu search \cite{SKORIN-KAPOV1990,TAILARD1991}, genetic method \cite{GAMBARDELLA1997} and 
simulated annealing \cite{CONNOLLY1990} 
have been developed. 
Those methods frequently attain near-optimal solutions, which 
often also happen to be the exact optimal solutions. The exactness 
is, however, not guaranteed in general. 

Most existing methods 
for finding the exact solutions of QAPs are designed with  
 the branch-and-bound (B\&B) method
\cite{ANSTREICHER2002,CLAUSEN1997,GONCALVES2015,HAHN2013,PARDALOS1997,ROUCAIROL1987}. 
As its name indicates, branching and (lower) bounding are the main procedures of the method. 
The bounding  based on doubly nonnegative (DNN) relaxations has recently attracted a great deal of attention as it provides tight bounds. 
Among the methods 
for solving large scale DNN relaxation problems,  
we mention four types 
of methods that can be applied to such DNN problems. 

First two of the methods 
are SDPNAL+~\cite{YST2015} 
(a majorized semismooth Newton-CG augmented Lagrangian method for semidefinite programming with nonnegative constraints) 
and BBCPOP~\cite{ITO2019}  
(a bisection and projection method for Binary, Box and Complementary constrained Polynomial Optimization Problems). 
Some numerical results on these two methods applied to QAP instances with dimensions $n=15$ to $50$ from 
QAPLIB~\cite{QAPLIB} were reported in \cite{ITO2019}. 
where BBCPOP attained tighter lower bounds for many of instances with dimensions $30$ to $50$ in 
less execution time. 
In addition, new and improved lower bounds were computed by Mittelmann using BBCPOP for the unknown minimum 
values of larger scale QAP instances including tai100a and tai100b. See QAPLIB~\cite{QAPLIB}. 
The third 
method 
is an alternating direction method of multipliers (ADMM) 
proposed by Oliveira et al.~\cite{OLIVEIRA2018} in combination with the facial reduction  for robustness. 
The fourth method is the Newton-bracketing (NB) method~\cite{ARIMA2018,KIM2019a}, 
which was developed to further improve 
the lower bounds  obtained from BBCPOP by incorporating the Newton method into BBCPOP. 


Numerical results on the NB method and BBCPOP are given in Table 2 
of~\cite{KIM2019a}, and ones on ADMM in Table 7 
of~\cite{OLIVEIRA2018}. Lower bounds for the QAP instances sko42, sko49, sko56, sko64, sko72, sko81, tai60a, 
and tai80a were commonly computed by these methods. 
While the lower bounds obtained by the NB method for the first three instances are not smaller than those by ADMM and the difference is as big as $1$, 
the NB method clearly attained the tightest lower bound among the three 
methods for the last five larger-scale instances with dimension $n \geq 60$. 

RLT~\cite{SHRARI1992} 
is known as a powerful method for global optimization of polynomial programming problems. 
In their paper \cite{GONCALVES2015}, Goncalves et al. implemented the level 2 RLT in their B\&B method 
run on heterogeneous CPUs and GPUs environment. They solved tai35b and tai40b from QAPLIB for the first time. 


The main motivation of our project is to
challenge larger scale QAP instances from QAPLIB~\cite{QAPLIB} that have not 
been solved yet. We implement the NB method 
\cite{KIM2019a} combined with the B\&B method in the specialized Ubiquity Generator (UG) framework~\cite{ShinanoAchterbergBertholdHeinzKoch2012} 
to find the exact solutions of large scale QAP instances.




The NB method has the following advantages: 
\vspace{-2mm}
\begin{description}
\item{(a) } The NB method generates a sequence of intervals $[\lb^k,\ub^k]$ 
$(k=1,2,\ldots)$ such that each $[\lb^k,\ub^k]$ contains the 
lower bound $\underline{\zeta}$ to be computed and both $\lb^k$ and $\ub^k$ converge monotonically to $\underline{\zeta}$ as 
$k \rightarrow \infty$. Hence, if $\lb^k$ becomes larger than the incumbent objective value $\hat{\zeta}$ at some iteration $k$, 
then $\hat{\zeta} < \lb^k \leq \underline{\zeta}$ follows.
Consequently, the node can be pruned. If $\ub^k < \hat{\zeta}$ occurs at some iteration $k$, then 
branching  the node should be determined as 
$\underline{\zeta} \leq \ub^k < \hat{\zeta}$ holds. In these two cases, the 
NB method 
can  terminate at iteration $k$ before the interval attains 
an accurate lower bound. This feature of the NB method 
works effectively to reduce a considerable amount of the computational time. 
\vspace{-2mm}
\item{(b) } The lower bounds computed by the NB method (and also the lower bounds computed by BBCPOP)
  are reliable or valid 
in the sense that they are guaranteed by very accurate dual feasible solutions of the Lagrangian DNN relaxation problem 
(see~\cite[Section 2.3]{KIM2019a}.) It would be plausible to  incorporate 
the technique used there for ADMM's lower bound, but 
the resulting lower bound would deteriorate. \vspace{-2mm}
\end{description}

The UG is a generic framework to
parallelize branch-and-bound based solvers and has achieved
large-scale MPI parallelism with 80,000 cores~\cite{SHINANO2016}.
It has also been generalized to handle non-branch-and-bound based solvers. 
Its experimental implementation 
for solving Shortest Vector Problem (SVP)  handles
more than 100,000 cores~\cite{TATEIWA2020} .

\medskip

The aim of this article is to present fundamentals of our
joint project for solving larger scale QAPs,  
report the progress, 
and discuss future plans to further develop the methods. We mention that  
tai35b, tai40b, tai30a and sko42 from QAPLIB~\cite{QAPLIB} have been successfully solved by our method, 
where the last two instances could be solved for 
the first time. 
QAPs with dimension at least 50 are our next target problems.

In Section~\ref{sect:subQAP}, we first
reformulate QAP~\eqref{eq:QAP0} to QAP~\eqref{eq:QAP}, which is  
a binary quadratic optimization problem form 
in $(1+n^2)$-dimensional vector variable $\u$. Then we 
introduce a class of its sub QAPs and 
describe the enumeration tree of those subproblems used in the branching process. 
In Section~\ref{sect:lowerBound}, a simple Lagrangian DNN relaxation problem 
of a sub QAP and the NB method to solve the relaxation problem are presented for 
the lower bounding procedure. In Section~\ref{sect:upperBound}, a method to retrieve a 
feasible solution of~\eqref{eq:QAP} from an approximate optimal solution of the Lagrangian DNN 
relaxation problem is described for the upper bounding procedure. 
We present three types branching rules in Section~\ref{sect:branchingRule}, and preliminary numerical 
results in Section~\ref{sect:numericalResults}. Some additional 
techniques are proposed for future development 
in Section~\ref{sect:additionalTechniques}.

\subsection*{Notation and Symbols}

Let $\Real^m$ denote the $m$-dimensional Euclidean space of column vectors $\z  = [z_1;\ldots;z_m]$, 
$\Real^m_+$ the nonnegative orthant of $\Real^m$, $\SymMat^m$ the linear space of $m\times m$ 
symmetric matrices with the inner product $\inprod{\A}{\B}  = \sum_{i=1}^m \sum_{j=1}^mA_{ij} B_{ij}$, 
and $\SymMat^m_+$ the cone of positive semidefinite matrices in $\SymMat^m$. 
For every $\ell \times m$ matrix $\U$, 
$\U._{j}$ denotes the $j$th column vector of $\U$, and vec~$\U$ the $\ell m$-dimensional column vector 
converted from the matrix $\U$, 
{\it i.e.}, 
\vec $\U = [\U._1;\ldots;\U._m]$. 
$\U^T$ denotes the transposition of $\U$. 
Throughout the paper, we use $n$ for the size of QAP~\eqref{eq:QAP0}, and $N = \{1,\ldots,n\}$ for the 
set of facilities and the set of locations. 
When $\Real^{\{0,N \times N\}}$ is used, each $(1+n^2)$-dimensional column vector 
$\u \in \Real^{\{0,N \times N\}}$ is 
represented component-wisely as 
$[u_0;u_{11},\ldots;u_{n1};u_{12};\ldots;u_{n2};\ldots,u_{1n};\ldots,u_{nn}]$, 
where $\U = [u_{ik}]$ forms an $n \times n$ matrix and 
$\u = [u_0;\mbox{vec }\U] = [u_0;\U._1;\ldots;\U._n]$. 
Each $((1+n^2) \times (1+n^2))$-dimensional symmetric 
matrix $\X \in \SymMat^{\{0,N \times N\}}$ has elements 
$X_{\salpha\sbeta}$ $(\balpha,\bbeta \in \{0,N \times N\})$.  


%

\section{A class of sub QAPs and the enumeration tree}

\label{sect:subQAP}

\subsection{Conversion of QAP~\eqref{eq:QAP0} into a QOP with a binary vector variable}

Let $\U = [u_{ij}]$ be an $n \times n$ matrix variable and 
$\u = [u_0;\vec(\U)] = [u_0;\U._{1};\ldots;\U._n] \in \Real^{\{0,N\times N\}} $ 
a $(1+n^2)$-dimensional column vector variable, where  
\begin{eqnarray*}
u_{ij} = \left\{\begin{array}{ll}
1 & \mbox{if facility $i$ is assigned to location $j$}, \\
0 & \mbox{otherwise}. 
\end{array}
\right.
\end{eqnarray*}
Define a constant matrix 
\begin{eqnarray*}
\Q^0 & = &
\begin{pmatrix}0 & \vec(\C)^T/2 \\ \vec(\C)/2 & \B\otimes \A \end{pmatrix}  \in 
\SymMat^{\{0,N \times N\}}, 
\end{eqnarray*}
where $\otimes$ denotes the Kronecker product. Then, the QAP can be expressed as 
\begin{eqnarray}
\zeta^0 & = & \min \left\{ 
\inprod{\Q^0}{\u\u^T}  : 
\begin{array}{l} 
\u \in \Real^{\{0,N\times N\}}_+, \  u_0 = 1,  \\ [3pt]
\displaystyle u_{ij}(u_0-u_{ij}) = 0 \ ((i,j) \in N\times N),
\\ [3pt]
\displaystyle  \sum_{k \in N} u_{ik} = u_0 \ (i \in  N), \ 
\displaystyle \sum_{k \in N} u_{kj} = u_0 \ (j\in N) 
\end{array}
 \right\}.  \label{eq:QAP} 
\end{eqnarray}

\subsection{A class of subproblems of QAP~\eqref{eq:QAP}}

\label{sec:subproblems}

We describe subproblems of QAP~\eqref{eq:QAP} with a family of 
subsequences $(i_1,\ldots,i_r)$ with $1 \leq i_1 < \cdots < i_r \leq n$. 
More precisely, let $\Pi_r$ denote the family of 
subsequences of $r$ distinct elements of $N$ $(r=1,\ldots,n)$. We assume $\Pi_0 = \{\emptyset\}$.  
$\Pi_n$ corresponds to the family of permutations of 
$1,\ldots,n$, and each $F =  (i_1,\ldots,i_r) \in \Pi_r$ 
(or $L =  (j_1,\ldots,j_r) \in \Pi_r$) with $r \geq 1$ corresponds to a permutation of 
 $r$ distinct elements $i_p$ $(p=1,\ldots,r)$ (or $j_p$ 
$(p=1,\ldots,r)$) of $N$.
For simplicity of notation, 
$F \in \Pi_r$ is frequently regarded as a subset of $N$ when 
 $i \in F$ and its complement 
$F^c = N\backslash F$ are described with respect to $N$. 

For each $(F,L) = ((i_1,\ldots,i_r), (j_1,\ldots,j_r)) \in \Pi_r \times \Pi_r$ $(r \in \{0,1,\ldots,n\})$,  let 
\begin{eqnarray*}
\Real^{\{0,F^c \times L^c\}} &=& \mbox{$1+|F^c||L^c|$-dimensional Euclidean space of column vectors $\x$} \\
& \ & \mbox{with elements $x_0$ and $x_{ij}$ $((i,j) \in F^c \times L^c)$}, \\
\Real^{\{0,F^c \times L^c\}}_+ & = & \mbox{the nonnegative orthant of $\Real^{\{0,F^c \times L^c\}}$}, \\ 
\SymMat^{\{0,F^c \times L^c\}} & = &  \mbox{the linear space of 
$(1+|F^c||L^c|) \times (1+|F^c||L^c|)$ symmetric} \\ 
& \ & \mbox{matrices $\X$ with elements $X_{\salpha\sbeta}$ $(\balpha,\bbeta \in \{0,F^c \times L^c\})$}, \\
\SymMat^{\{0,F^c \times L^c\}}_+ & = &  \mbox{the cone of positive semidefinite matrices in }
\SymMat^{\{0,F^c \times L^c\}}. 
\end{eqnarray*}
Then, for each $(F,L) = ((i_1,\ldots,i_r), (j_1,\ldots,j_r)) \in \Pi_r \times \Pi_r$ 
$(r \in \{0,1,\ldots,n\})$,
a subproblem of QAP~\eqref{eq:QAP} for assigning each facility $i_p \in F$
to each location $j_p \in L$ $(p=1,\ldots,r)$ 
can be written as 
\begin{eqnarray*}
\mbox{QAP$(F,L)$: }
\zeta(F,L) =
\min\left\{\inprod{\Q(F,L)}{\x\x^T}: 
\begin{array}{l}
\x \in \Real^{\{0,F^c \times L^c\}}_+, \ x_0 = 1, \\ [2pt] 
x_{ij}(x_0-x_{ij})=0   \\ [2pt]
((i ,j) \in F^c \times  L^c) \\ [2pt]
\displaystyle \sum_{k\in L^c}x_{ik}- x_0 = 0 \ (i \in F^c), \\ [2pt]
\displaystyle \sum_{k\in F^c}x_{kj}- x_0 = 0 \ (j\in L^c)
\end{array}
\right\} 
\end{eqnarray*}
for some matrix $\Q(F,L) \in \SymMat^{\{0,F^c \times L^c\}}$, which is given by $\Q(F,L) =  \P(F,L)^T \Q^0 \P(F,L)$,
where 
\begin{eqnarray*}
\P(F,L) & = & \mbox{the $(1+|N \times N|) \times (1+|F^c \times L^c|)$ matrix with elements}\\
& \ & P(F,L)_{\salpha\sbeta} \ = \ \left\{\begin{array}{ll}
1 & \mbox{if } \balpha = 0 \ \mbox{and } \bbeta = 0, \\ 
1 & \mbox{if } \balpha = (i_p, j_p) \text{ for some } p\in\{1, \ldots, r\} \text{ and } \bbeta=0, \\ 
1 & \mbox{if } \balpha = \bbeta \in F^c \times L^c, \\
0 & \mbox{otherwise. }
\end{array}\right.
\end{eqnarray*}
Obviously, $\P(\emptyset,\emptyset)$ coincides with the 
identity matrix in $\SymMat^{\{1,N \times N\}}$ and 
QAP$(\emptyset,\emptyset)$ corresponds to QAP~\eqref{eq:QAP}. 

\subsection{Enumeration tree} 

\label{sect:enumerationTree} 

A fundamental step in the B\&B method is to successively 
construct an enumeration tree. Each node of the enumeration tree 
corresponds to a subproblem of QAP~\eqref{eq:QAP}, QAP$(F,L)$ 
for some $(F,L) \in \Pi_r\times \Pi_r$ with $r \in \{0,1,\ldots,n\}$. Thus, each node is identified with a subproblem. 
 To initialize the tree, we set QAP$(\emptyset,\emptyset)$ 
({\it i.e.,} QAP~\eqref{eq:QAP})  for the root node and compute the incumbent objective value 
$\hat{\zeta}$ by applying a heuristic method to QAP~\eqref{eq:QAP}. 

Suppose that we have generated a node QAP$(F,L)$, 
where $F = (i_1,\ldots,i_r) \in \Pi_r$ and 
$L = (j_1,\ldots,j_r) \in \Pi_r$ for some $r\in \{1,\ldots,n\}$ 
(or $(F,L) = (\emptyset,\emptyset)$). 
Then, we compute a lower bound $\zeta^{\ell}(F,L)$ (see Section~\ref{sect:lowerBound}) 
and an upper bound $\zeta^{u}(F,L)$ (see Section~\ref{sect:upperBound}), 
which is associated with a feasible solution of QAP$(F,L)$, 
for  the (unknown) optimal value $\zeta(F,L)$; $\zeta^{\ell}(F,L) \leq \zeta(F,L) \leq \zeta^{u}(F,L)$.   
The incumbent objective value $\hat{\zeta}$ is then updated by 
$\hat{\zeta} = \min \{\hat{\zeta}, \ \zeta^{u}(F,L) \}$.  

When $\hat{\zeta} \leq \zeta^{\ell}(F,L)$ occurs, 
the (unknown) optimal value $\zeta(F,L)$ of QAP$(F,L)$ 
 is greater than or equal to 
the incumbent objective value $\hat{\zeta}$. In  this case, the node is pruned (or terminated). 
Otherwise, branching the node into $n-r$ child nodes takes place  as follows. 
We select either facility $f_+ \in F^c$ or 
location $\ell_+ \in L^c$. 
If  $f_+ \in F^c$ is selected, then 
$|L^c| = n-r$ child nodes $\mbox{QAP}((F,f_+),(L,\ell)) \ (\ell \in L^c)$ 
 are generated. Otherwise,   
$|F^c| = n-r$ child nodes 
$\mbox{QAP}((F,f),(L,\ell+)) \  (f \in F^c)$ are 
generated.  A branching rule determines which $f_+ \in F^c$ or $\ell_+ \in L^c$ to select.
In Section~\ref{sect:branchingRule}, we describe three types of branching rules in detail. 

For any QAP$(F,L)$ 
possibly located at the bottom level (or the $n$th level), 
we have $|F|=|L|=n$ and $F^c = L^c = \emptyset$. Thus, 
QAP$(F,L)$ becomes a 
trivial QAP with every facility assigned to some location and 
vice versa. It is trivial to compute its optimal value $\zeta(F,L)$. 
We note that the height of the enumeration tree is at most $n$.
Numerically, we can enumerate all $(n-r)!$ feasible solutions of QAP$(F,L)$ to compute its exact optimal value 
when $|F^c|=|L^c|=n-r \leq m$ holds for a small enough $m$. In the numerical results reported in Section~\ref{sect:numericalResults}, 
we took $m=7$.

\section{Lower bounding procedure}

\label{sect:lowerBound}

Throughout this section, we fix $(F,L) \in \Pi_r \times \Pi_r$ $(r\in \{0,1,\ldots,n\})$, and 
present a method to compute a lower bound for the optimal objective value $\zeta(F,L)$ 
of QAP$(F,L)$. 
Before deriving a doubly nonnegative (DNN) relaxation of QAP$(F,L)$ 
 in Section~\ref{sect:DNN} 
and a Lagrangian DNN relaxation of QAP$(F,L)$ 
in Section~\ref{sect:LagDNN}, we first 
 transform  QAP$(F,L)$ into $\mbox{\bQAP}(F,L)$ described below.

We note that each homogeneous linear equality 
constraint $\displaystyle \sum_{k\in L^c}x_{ik}- x_0 = 0$ of QAP$(F,L)$ 
can be written  as $\hat{\cc}(i,L)^T \x = 0$ for some $\hat{\cc}(i,L) \in \Real^{\{0,F^c\times L^c\}}$  $(i \in F^c)$. 
Let $\hat{\a}(i,L) = \vec(\hat{\cc}(i,L)\hat{\cc}(i,L)^T)^T$. Then, the linear equality 
is equivalent to $\hat{\a}(i,L) \vec(\x\x^T)  = 0$. 
Similarly, we can rewrite each homogeneous linear constraint 
$\displaystyle \sum_{k\in F^c}x_{kj}- x_0 = 0$ as 
$\tilde{\a}(j,F)\vec(\x\x^T) = 0$ for some $(1+|F^c||L^c|)^2$-dimensional row 
vector $\tilde{\a}(j,F)$ $(j \in L^c)$. Hence, the entire linear constraints 
\begin{eqnarray*}
& & \sum_{k\in L^c}x_{ik}- x_0 = 0 \ (i \in F^c) \ \mbox{ and } 
\sum_{k\in F^c}x_{kj}- x_0 = 0 \ (j \in L^c)
\end{eqnarray*}
can be expressed as $\A(F,L) \vec(\x\x^T) = \0 \in \Real^{|F^c| + |L^c|}$, 
where $\A(F,L)$ denotes the $(|F^c| + |L^c|) \times 
(1+|F^c||L^c|)^2$ matrix consisting of the row vectors 
$\hat{\a}(i,L)$ $(i \in F^c)$ and $\tilde{\a}(j,F)$ $(j \in L^c)$. 

Define 
\begin{eqnarray*}
\coneK_1(F,L) & = & \SymMat^{\{0,F^c\times L^c\}}_+,
\\ [2pt]
\coneK_2(F,L) & = & \left\{ \X \in \SymMat^{\{0,F^c\times L^c\}}  : 
\begin{array}{l}
X_{\salpha\sbeta} \geq 0 \ (\balpha,\bbeta \in \{0,F^c\times L^c\}), \\ [2pt]
X_{0\salpha}  - X_{\salpha\salpha} = 0 \ (\balpha \in F^c\times L^c)
\end{array} 
\right\}, \\ 
\H(F,L) & = & 
\mbox{the matrix in $\SymMat^{\{0,F^c\times L^c\}}$ with $1$ at the $(0,0)$th  element}\\ 
& \ & \mbox{and $0$ elsewhere}. 
\end{eqnarray*}
By construction, 
the constraint $x_{ij}(x_0 - x_{ij}) = 0$ $((i,j) \in F^c \times L^c))$ is equivalent 
to $X_{0\salpha}  - X_{\salpha\salpha} = 0 \ (\balpha \in F^c\times L^c)$ 
if $\X = \x\x^T \in \SymMat^{\{0,F^c\times L^c\}}$ and $\x \in \Real^{\{0,F^c\times L^c\}}$, and 
$\inprod{\H(F,L)}{\X} = X_{00}$ for every $\X \in \SymMat^{\{0,F \times L\}}$.  
Therefore, we can rewrite QAP$(F,L)$ given in Section~\ref{sec:subproblems} as
\begin{eqnarray*}
\mbox{\bQAP}(F,L):  
\zeta(F,L)  \hspace{-3mm} & = & \hspace{-3mm}  \left\{\inprod{ \Q(F,L)}{\X} 
\hspace{-1mm} : \hspace{-2mm} 
\begin{array}{l}
\x \in \Real^{\{0,F^c\times L^c\}}, \ 
\X = \x\x^T \in \coneK_1(F,L)\cap\coneK_2(F,L),\\[2pt]
\inprod{\H(F,L)}{\X}  = 1, \ 
\A(F,L) \vec(\X) = \0
\end{array} 
\right\}.
\end{eqnarray*}

\subsection{A DNN relaxation of \bQAP$(F,L)$} 

\label{sect:DNN} 

In \bQAP$(F,L)$, the constraint $\X = \x\x^T$ requires 
the rank of the variable matrix $\X$ to be $1$. By removing this constraint, we 
obtain 
\begin{eqnarray*}
\mbox{DNN}^p(F,L):  
\eta^p(F,L)  & = & \min 
 \left\{\inprod{ \Q(F,L)}{\X} 
\hspace{-1mm} :  
\begin{array}{l}
\X \in \coneK_1(F,L)\cap\coneK_2(F,L),\\[2pt]
 \inprod{\H(F,L)}{\X} = 1, \ 
\A(F,L) \vec(\X) = \0
\end{array} 
\right\}.
\end{eqnarray*}
which serves as a DNN relaxation of \bQAP$(F,L)$
with $\eta^p(F,L) \leq \zeta(F,L)$. 

\subsection{A  Lagrangian DNN relaxation problem of \bQAP$(F,L)$}

\label{sect:LagDNN}

We know that 
\begin{eqnarray*}
\A(F,L) \vec(\X) \geq \0 \ \mbox{for every } \X \in \coneK_1(F,L) = \SymMat^{\{0,F^c\times L^c\}}_+, 
\label{eq:nonnegative}
\end{eqnarray*}
since each row vector of $\A(F,L)$ is of the form $\vec(\cc\cc^T)$ for some 
$\cc \in \Real^{\{0,F^c\times L^c\}}$.  
Hence, the constraint $\A(F,L) \vec(\X) = \0$ in DNN$^p(F,L)$  can be replaced by 
\begin{eqnarray*}
0 & = & \e^T \A(F,L) \vec(\X) \ = \ \Inprod{\sum_{i\in F^c} \hat{\cc}(i,L)  \hat{\cc}(i,L)^T + 
\sum_{j\in L^c} \tilde{\cc}(j,F)  \tilde{\cc}(i,L)^T }{\X },  
\end{eqnarray*}
where $\e \in \Real^{|F|+|L|}$ denotes the vector of $1$'s,  and the second equality holds 
by the construction of $ \A(F,L)$. 
Applying the Lagrangian relaxation to the resulting problem, 
we obtain 
\begin{eqnarray*}
\mbox{DNN}^p_{\lambda}(F,L):  
\eta^p_{\lambda}(F,L)  & = & \min \left\{\inprod{ \Q_{\lambda}(F,L)}{\X}  : 
\begin{array}{l}
\X \in \coneK_1(F,L)\cap\coneK_2(F,L), \\ 
 \inprod{\H(F,L)}{\X} = 1 
\end{array} 
\right\}, 
\end{eqnarray*}
and its dual
\begin{eqnarray*}
\mbox{DNN}^d_{\lambda}(F,L):  
\eta^d_{\lambda}(F,L)  & = & \max\left\{ y \in \Real : 
\begin{array}{l}
\Y_1 \in \coneK_1(F,L)^*, \ \Y_2 \in \coneK_2(F,L)^*,\\[2pt]
\Q_{\lambda}(F,L) - \H(F,L)y = \Y_1 + \Y_2
\end{array}
\right\}. 
\end{eqnarray*}
Here $\lambda \in \Real$ denotes a Lagrangian multiplier, and 
\begin{eqnarray*}
\Q_{\lambda}(F,L) & = & \Q(F,L) + \lambda
\left(\sum_{i\in F^c} \hat{\cc}(i,L)  \hat{\cc}(i,L)^T + 
\sum_{j\in L^c} \tilde{\cc}(j,F)  \tilde{\cc}(i,L)^T \right) \in \SymMat^{\{0,F^c\times L^c\}}.
\end{eqnarray*} 
DNN$^p_{\lambda}(F,L)$ serves as a 
Lagrangian relaxation of DNN$^p(F,L)$ 
and a Lagrangian DNN relaxation of \bQAP$(F,L)$ with 
 $\eta^p_{\lambda}(F,L) \leq \eta^p(F,L) \leq \zeta(F,L)$ for any $\lambda \in \Real$. 
 We can prove that $\eta^p_{\lambda}(F,L)$ converges monotonically to 
 $\eta^p(F,L)$ as $\lambda \rightarrow \infty$. 
 We also see by the duality theorem 
 (\cite[Lemma 2.3]{ARIMA2018}) that $\eta^p_{\lambda}(F,L) = \eta^d_{\lambda}(F,L)$. 

\subsection{The Newton-bracketing method for solving the 
primal-dual pair of DNN$^p_{\lambda}(F,L)$ and DNN$^d_{\lambda}(F,L)$} 

\label{sec:NewtonBracket}

To describe the NB method for solving DNN$^p_{\lambda}(F,L)$ 
and DNN$^d_{\lambda}(F,L)$, $\lambda$ is fixed to a sufficiently large positive number, say 
$\lambda = 100,000$. 
The optimal value $\eta^d_{\lambda}$ of DNN$^d_{\lambda}(F,L)$ 
is denoted by $y^*$ throughout this subsection.

Let $-\infty = \lb^0 < y^* < \ub^0$. 
The NB method  applied to DNN$^p_{\lambda}(F,L)$ 
and DNN$^d_{\lambda}(F,L)$ generates 
\begin{eqnarray*}
& & \widehat{\X}^k \in \coneK_1(F,L)\cap\coneK_2(F,L), \ 
(y^k,\widehat{\Y}_1^k, \widehat{\Y}_2^k) 
\in \Real \times \SymMat^{\{0,F^c \times L^c\}} \times \coneK_2(F,L)^*,  \\[3pt]
& & \lb^k \in \Real, \  \ub^k \in \Real, \ 0 \geq \tau^k \in \Real 
 \end{eqnarray*}
$(k=1,2,\ldots)$ such that 
\begin{eqnarray}
\left.
\begin{array}{l}
\Q(F,L) - \H(F,L) y^k - \widehat{\Y}^k_1 -  \widehat{\Y}^k_2 = \O, \\[3pt]
\lb^k = y^k + (1+|F^c|)\tau^k, \ 
\tau^k = \min\{0,\mbox{the minimum eigenvalue of $\widehat{\Y}^k_1$}\}, \\[3pt]
 \lb^0 \leq \lb^k \leq \lb^{k+1} \rightarrow y^* \leftarrow \ub^{k+1} \leq \ub^k \leq \ub^0, \\[3pt] 
\X = \widehat{\X}^k \ \mbox{is a feasible solution of~DNN$^p_{\lambda}(F,L)$},    \ 
\inprod{\Q(F,L)}{\widehat{\X}^k} = \ub^k. 
\end{array}
\right\} \label{eq:NBiterate}
\end{eqnarray}
We see from these relations that for sufficiently large $k$, $\X = \widehat{\X}^k$ 
and $(y,\Y_1,\Y_2) = (y^k,\widehat{\Y}_1^k,\widehat{\Y}_2^k)$ 
provides approximate optimal solutions 
of~DNN$^p_{\lambda}(F,L)$ and DNN$^d_{\lambda}(F,L)$, respectively. Hence 
$\X = \widehat{\X}^k$  is  an approximate optimal solution of  ~DNN$^p(F,L)$. 
See \cite{KIM2019a} for more details. 

It is natural to use a  stopping criterion for the NB method as 
 \vspace{-2mm}
\begin{description}
\item{(a) } $|\ub^k-\lb^k| < \epsilon \max\{|\lb^k|,|\ub^k|,1\}$, where $\epsilon$ denotes a sufficiently small positive number. 
 \vspace{-2mm}
\end{description}
When  the NB method  is implemented in the B\&B method for QAPs with integer optimal values, 
an accurate 
approximation of $y^*$  does not  have to be computed all the time by the NB method. 
We can terminate the iteration in case when 
\vspace{-2mm}
\begin{description}
\item{(b) } $\hat{\zeta} \leq \lceil\lb^k\rceil$ or \vspace{-2mm}
\item{(c) } $\ub^k < \hat{\zeta}$, 
\vspace{-2mm}
\end{description}
where $\hat{\zeta}$ denotes 
an integer incumbent objective  value of the root node QAP. 
If (b) occurs, then we know 
that the sub QAP, \bQAP$(F,L)$ to which the NB method has applied, 
does not have any feasible solution whose 
objective value  is smaller than the incumbent objective value 
$\hat{\zeta}$; hence \bQAP$(F,L)$ can be pruned. If case (c) occurs, 
then $y^* \leq \ub^k < \hat{\zeta}$. 
Thus, branching  \bQAP$(F,L)$ further is necessary even if  the NB 
iteration continues. By employing these two additional criteria (b) and (c), we can save a lot of computational time for the 
NB method. This is a distinctive advantage of using the NB method in the B\&B method.

\section{Upper bounding procedures}

\label{sect:upperBound}

Let $(F,L) \in \Pi_r \times \Pi_r $ for some $r \in \{0,1,\ldots,n\}$. 
An approximate optimal solution of QAP$(F,L)$ needs to be 
computed for the upper 
bounding procedure. 
A simple rounding of an approximate optimal 
solution of its Lagrangian relaxation DNN$^p_{\lambda}(F,L)$ and the tabu search 
method~\cite{SKORIN-KAPOV1990,TAILARD1991} for QAP$(F,L)$  have been used for the computation. 
We plan to replace this method by the one presented below in this section to  effectively 
 utilize information from DNN$^p_{\lambda}(F,L)$. 
In the remainder of this section, we assume  $F=L=\emptyset$ 
so that QAP$(F,L)$ coincides with the root node QAP~\eqref{eq:QAP}, for simplicity of notation.
All the discussions, however, can be easily adapted to any sub QAP, QAP$(F,L)$.

Suppose that we have applied the NB method to the primal-dual pair of DNN problems 
DNN$^p_{\lambda}(\emptyset,\emptyset)$ and DNN$^d_{\lambda}(\emptyset,\emptyset)$ and obtained  
an approximate optimal solution $\widehat{\X}$ of DNN$^p_{\lambda}(\emptyset,\emptyset)$.
The set of feasible solutions $\X$ of 
DNN$^p_{\lambda}(\emptyset,\emptyset)$ may be 
regarded as a doubly nonnegative relaxation of  
the set of $\u\u^T$ with feasible solutions $\u$  of 
QAP~\eqref{eq:QAP}, and the first column 
of $\widehat{\X}$ corresponds to a feasible (or approximate optimal) 
solution $u_0\u = \u$ of QAP~\eqref{eq:QAP}. We know that 
$\U = [u_{ik}]$ forms a permutation matrix for every feasible 
solution $\u$ of QAP~\eqref{eq:QAP}. Through preliminary 
numerical experiment, we have observed that the $n \times n$ 
matrix
\begin{eqnarray*}
\widehat{\U} & = & 
\begin{pmatrix} 
\widehat{X}_{11,0} & \widehat{X}_{12,0} & \ldots & \widehat{X}_{1n,0} \\
\widehat{X}_{21,0} & \widehat{X}_{22,0} & \ldots& \widehat{X}_{2n,0} \\
\cdots & \cdots & \cdots & \cdots \\
\widehat{X}_{n1,0} & \widehat{X}_{n2,0} & \ldots & \widehat{X}_{nn,0} 
\end{pmatrix}
\end{eqnarray*}
approximately forms a doubly-stochastic matrix. 
A popular approach for recovering a feasible solution 
from 
$\widehat{\U}$ is through  appropriate rounding procedures. 
Here we propose a method 
to compute a permutation matrix $\U = \overline{\U}$ which minimize the 
distance $\parallel\U - \widehat{\U}\parallel$ among all permutation matrices $\U$: 
\begin{eqnarray*}
\bar{\zeta} & = & \min \left\{ \parallel \U -\widehat{\U}\parallel^2 : \mbox{ $\U \in \Real^{n \times n}$ is a permutation matrix }  \right\} \\
& = & \min \left\{  -2\inprod{\widehat{\U}}{\U}  + \parallel\widehat{\U}\parallel^2 
+ n^2 : \mbox{$\U \in \Real^{n \times n}$ is a permutation matrix }  \right\} \\
& \ & \mbox{(since $\parallel\U\parallel^2 = n^2$ for every permutation 
matrix $\U\in \Real^{n \times n}$)}.
\\
& = & \min \left\{  \inprod{-2\widehat{\U}}{\U} : \mbox{$\U \in \Real^{n \times n}$ is a permutation matrix }  \right\}   + \parallel\widehat{\U}\parallel^2 + n^2,  
\end{eqnarray*} 
where $\parallel\cdot\parallel$ denotes the Frobenius norm. 
This problem forms a linear assignment problem, which can be solved the well-known 
Hungarian method \cite{KUHN1955}.

\section{Branching rules}

\label{sect:branchingRule}

Let $(F,L) \in \Pi_r \times \Pi_r$ for some $r \in \{0,1,\ldots,n\}$. 
As we have briefly mentioned in Section~\ref{sect:enumerationTree},
if node QAP$(F,L)$ is not pruned, we select either an $f_+ \in F^c$ to branch 
QAP$(F,L)$ into $|L^c|$ child nodes QAP$((F,f_+),(L,\ell))$ $(\ell \in L^c)$ or an $\ell_+ \in L^c$ to branch 
QAP$(F,L)$ into $|F^c|$ child nodes QAP$((F,f),(L,\ell_+))$ $(f \in F^c)$. 

We present three types of rules 
for selecting an $f_+ \in F^c$ or an $\ell_+ \in L^c$. In each rule, an 
evaluation function $\varphi(f,\ell)$ is introduced for possible child nodes QAP$((F,f),(L,\ell))$ 
$((f,\ell) \in F^c \times L^c)$ so that $\varphi(f,\ell)$ can represent 
the optimal values of 
QAP$((F,f),(L,\ell))$ or its Lagrangian DNN relaxation problem DNN$^p_{\lambda}((F,f),(L,\ell))$.
Then, the following functions are  defined:
\begin{eqnarray*}
\bar{\varphi}(f,\cdot) & = & \frac{1}{|L^c|}\sum_{\ell \in L^c} \varphi(f,\ell) \ \
\mbox{(the mean of $\varphi(f,\ell)$ over $\ell \in L^c$) } (f \in F^c), \\ 
\bar{\varphi}(\cdot,\ell) & = & \frac{1}{|F^c|}\sum_{f \in F^c} \varphi(f,\ell) \ \
\mbox{(the mean of $\varphi(f,\ell)$ over $f \in F^c$) } (\ell \in L^c). 
\end{eqnarray*}
Our method is to choose $f_+ = \mbox{arg max} \{\bar{\varphi}(f,\cdot): f \in F^c\}$ 
if $\mbox{max} \{\bar{\varphi}(f,\cdot): f \in F^c\} 
\geq \mbox{max} \{\bar{\varphi}(\cdot,\ell): \ell \in L^c\}$, 
and $\ell_+ = \mbox{arg max} \{\bar{\varphi}(\cdot,\ell): \ell \in L^c\}$ otherwise. 
Each rule employs a different function $\varphi(f,\ell)$ from others. We describe 
how $\varphi(f,\ell)$  is defined with a fixed $((f,\ell) \in F^c \times L^c)$ for the three rules 
in Sections~\ref{sect:QAPmean}, \ref{sect:DNNprimal} 
and~\ref{sect:DNNprimal}, respectively.

\subsection{Branching rule M: a rule using the mean of all feasible objective values of 
Q$(F,L)$ $((F,L) \in \Pi_r \times \Pi_r, \ r=0,1,\ldots,n)$} 

\label{sect:QAPmean}

Define $\varphi(f,\ell)$ to be 
the mean of the objective values of QAP$((F,f),(L,\ell))$ over its feasible solutions, which can be  
computed by substituting $x_0 = 1$ and $x_{ij} = 1/(n-r-1)$  
$((i,j) \in (F,f)^c \times (L,\ell)^c)$ into the objective function 
$\inprod{\Q((F,f),(L,\ell))}{\x\x^T}$ of QAP$((F,f),(L,\ell))$. 
This branching rule M is employed in the current parallel implementation of the \BB \ 
method on the UG framework. 

\subsection{Branching rule P: 
a rule using an approximate optimimal solution $\widehat{\X}$ of 
DNN$^p_{\lambda}(F,L)$}

\label{sect:DNNprimal}

Let $\widehat{\X}$ be an approximate optimal solution of DNN$^p_{\lambda}(F,L)$.
Then, $\widehat{\X}$ may be regarded 
as an approximate optimal solution of DNN$^p(F,L)$ since 
DNN$^p_{\lambda}(F,L)$ is a Lagrangian relaxation of DNN$^p(F,L)$ and its optimal 
value $\eta^p_{\lambda}(F,L)$ converges to $\eta^p(F,L)$ of DNN$^p(F,L)$ as $\lambda \rightarrow \infty$. 
We will construct a rough approximate optimal solution $\widehat{\X}(f,\ell)$ 
of DNN$^p((F,f),(L,\ell))$ from $\widehat{\X}$, and define $\varphi(f,\ell)$ to be the objective value 
of DNN$^p((F,f),(L,\ell))$ at $\widehat{\X}(f,\ell)$, 
{\it i.e.}, $\varphi(f,\ell) = \inprod{\Q((F,f),(L,\ell))}{\widehat{\X}(f,\ell)}$ as described below. 

We first project $\widehat{\X} \in \SymMat^{\{0,F^c \times L^c\}}$ onto the linear space 
$\SymMat^{\{0,(F,f)^c \times (L,\ell)^c\}}$ where the feasible region of QAP$^p((F,f),(L,\ell))$ 
lies. Let ${\X}(f,\ell)$ denote the metric projection of $\widehat{\X}$, 
 which is obtained by deleting the rows $\widehat{\X}_{i.}$ 
$(i \not\in \{0,(F,f)^c \times (L,\ell)^c\})$ and the columns 
$\widehat{\X}_{.j}$ $(j \not\in \{0,(F,f)^c \times (L,\ell)^c\})$  from $\widehat{\X}$ simultaneously. 
We also know that every feasible solution $\X$ of DNN$^p((F,f),(L,\ell))$ must lie in the linear manifold 
\begin{eqnarray*}
M & = & \left\{ \X \in \SymMat^{\{0,(F,f)^c \times (L,\ell)^c\}} : 
\begin{array}{l}
X_{00} = 1, \ \A((F,f),(L,\ell)) \vec(\X) = \0, \\[3pt]
X_{0\salpha} - X_{\salpha\salpha} = 0 \ (\balpha \in (F,f)^c \times (L,\ell)^c)
\end{array}
\right\}.  
\end{eqnarray*}
Let $\widehat{\X}(f,\ell)$ be the metric projection of ${\X}(f,\ell)$ onto $M$, and then 
define 
\[
\varphi(f,\ell) = \inprod{\Q((F,f),(L,\ell))}{\widehat{\X}(f,\ell)}.
\] 
We note that $\widehat{\X}(f,\ell)$ is not necessarily in 
$\coneK_1((F,f),(L,\ell))  = \SymMat^{\{0,(F,f)^c \times (L,\ell)^c\}}_+$; hence $\widehat{\X}(f,\ell)$ 
is not necessarily a feasible solution of DNN$^p((F,f),(L,\ell))$ in general. 

\subsection{Branching rule D: a rule using an approximate optimal solution 
$(\hat{y},\widehat{\Y}_1,\widehat{\Y}_2)$ of 
DNN$^d_{\lambda}(F,L)$}

\label{sect:DNNdual}

Recall that the NB method applied to the primal-dual pair of DNN$^p_{\lambda}(F,L)$ 
and DNN$^d_{\lambda}(F,L)$ generates 
\begin{eqnarray*}
& & \X^k \in \coneK_1(F,L)\cap\coneK_2(F,L), \ 
(y^k,\Y_1^k,\Y_2^k) 
\in \Real \times \SymMat^{\{0,F \times L\}} \times \coneK_2(F,L)^*,  \\
& & \lb^k \in \Real, \  \ub^k \in \Real, \ 0 \geq \tau^k \in \Real 
\end{eqnarray*}
$(k=1,2,\ldots)$ 
satisfying~\eqref{eq:NBiterate}. We assume that the method terminates with the stopping criterion (c) 
$\ub^k < \hat{\zeta}$, which requires to branch the node QAP$(F,L)$. 

Define the 
$(1+|F^c \times L^c|) \times (1+|(F,f)^c \times (L,\ell)^c|)$ matrix 
$\P(f,\ell)$ by 
\begin{eqnarray*}
P(f,\ell)_{\salpha\sbeta} & = & \left\{\begin{array}{ll}
1 & \mbox{if } \balpha = 0 \ \mbox{and } \bbeta = 0, \\ 
1 & \mbox{if } \balpha = (f,\ell) \ \mbox{ and } \bbeta = 0, \\ 
1 & \mbox{if } \balpha = \bbeta \in (F,f)^c \times (L,\ell)^c, \\
0 & \mbox{otherwise }.
\end{array}\right.
\end{eqnarray*}
Then 
\begin{eqnarray*}
& & \Q((F,f),(L,\ell)) = \P(f,\ell)^T \Q(F,L) \P(f,\ell) \in \SymMat^{\{0,(F,f) \times (L,\ell)\}}, \\ 
& & \H((F,f),(L,\ell)) = \P(f,\ell)^T \H(F,L) \P(f,\ell) \in  \SymMat^{\{0,(F,f) \times (L,\ell)\}}, \\ 
& & \widehat{\Y}_1(f,\ell) \equiv  \P(f,\ell)^T\Y^k_1 \P(f,\ell) \in \SymMat^{\{0,(F,f)^c\times(L,\ell)^c\}}, \\ 
& & \widehat{\Y}_2(f,\ell) \equiv  \P(f,\ell)^T\Y^k_2 \P(f,\ell) \in \coneK_2(F,f),(L,\ell))^*
\end{eqnarray*}
hold. 
We apply the above to the identity $\Q(F,L) - \H(F,L) y^k - \Y^k_1 -  \Y^k_2 = \O$ 
in~\eqref{eq:NBiterate}. More precisely,
multiplying $\P((F,f),(L,\ell))^T$ to each term of the identity on the left 
and $\P((F,f),(L,\ell))$ on the right results in 
\begin{eqnarray*}
& & \O \ = \ \Q((F,f),(L,\ell)) - \H((F,f),(L,\ell))y^k - \widehat{\Y}_1(f,\ell) - \widehat{\Y}_2(f,\ell), \\[2mm]
& \ & \widehat{\Y}_1(f,\ell) \in \SymMat^{\{0,(F,f)^c\times(L,\ell)^c\}}, \ 
\widehat{\Y}_2(f,\ell) \in \coneK_2(F,f),(L,\ell))^*. 
\end{eqnarray*}
Now, define
\begin{eqnarray*}
\tau(f,\ell) & = &\min\{0,\mbox{the minimum eigenvalue of $\widehat{\Y}_1(f,\ell)$}\}, \\ 
\varphi(f,\ell) & = & y^k + (1+|(F,f)^c|)\tau(f,\ell) \ = \ y^k +(n-r)\tau(f,\ell). 
\end{eqnarray*}
Then we can prove that $\lb^k \leq \varphi(f,\ell) \leq \eta^p_{\lambda}((F,f),(L,\ell))$,  
which implies that $\varphi(f,\ell)$ itself serves as a lower bound of the optimal value 
$\zeta((F,f),(L,\ell))$ of the child node QAP$((F,f),(L,\ell))$. Thus, 
if $\hat{\zeta} \leq \varphi(f,\ell)$ holds, then we can terminate 
QAP$((F,f),(L,\ell))$ without applying the NB method to QAP$((F,f),(L,\ell))$. This is a 
distinct advantage of the branching rule {\bf D},  although the bound 
$\varphi(f,\ell)$ 
may not be as strong as $\eta^p_{\lambda}((F,f),(L,\ell))$.

\section{Preliminary numerical results}

\label{sect:numericalResults}

We have implemented two types of codes for the parallel \BB method 
to solve QAPs. 
The first one is written in MATLAB, which has been developed to see how effectively 
and efficiently the NB method works in the branch-and-bound framework. 
The code used there for the NB method 
is modifications of the ones written for the MATLAB 
software package NewtBracket \cite{KIM2020}, which was released very recently. 
Specifically, the stopping criteria of the NB method were changed as 
presented in the last paragraph of Section~\ref{sec:NewtonBracket}. 
All techniques presented in Sections~\ref{sect:subQAP}, \ref{sect:lowerBound} 
and~\ref{sect:branchingRule} have been incorporated in the MATLAB code, but not 
the upper bounding procedure described in Section~\ref{sect:upperBound}.  
The MATLAB code works in parallel but can solve only small-sized  QAP instances 
as shown in Section~\ref{sect:MATLAB}. 

The second code is a C++ implementation of the MATLAB code.
Based on the information and knowledge obtained from preliminary numerical 
experiment by the MATLAB code, we have been developing the C++ code running on the UG, a generic framework to parallelize branch-and-bound based 
solvers.
An application of the UG can be developed on PC, in which communication is carried out through shared memory~\cite{SHINANO2013}. 
After  recompiling the application code with additional small amount of MPI related code for transferring objects, the parallel solver could potentially run with 
more than 100,000 cores~\cite{SHINANO2016, TATEIWA2020}.

Two types of processes exist when running the parallelized QAP solver by UG on distributed memory environment.
First, there is a single \textsc{LoadCoordinator}, 
which makes all decisions
concerning the dynamic load balancing and distributes subproblems of
QAP instances to be solved.
Second, all other processes are \textsc{Solver}s 
that solve the distributed subproblem by regarding 
it as root node. 
The UG tries to solve all open subproblems on the single branch-and-bound search tree in parallel as 
much as it can by using all available \textsc{Solver}{s}. 
More precisely, 
the UG executes dynamic load balancing so that subproblems of the branch-and-bound search tree are transferred to the other \textsc{Solver}{s} 
when an idle \textsc{Solver}\ 
exists.
In order to reduce the idle time of the \textsc{Solver}{s}, 
\textsc{LoadCoordinator}\ 
tries to keep a small number 
of open subproblems 
(the initial number can be specified by a runtime parameter and it is changed dynamically during the computation) so that it can assign a subproblem at the earliest time 
when an idle \textsc{Solver}\ 
exists. 

For {\em ramp-up},  a phase until all cores become busy (see \cite{RALPHS2018}),
 several methods exist. 
In our case of the parallelized QAP solver, we used the {\em normal ramp-up}, which performs the 
normal parallel branch-and-bound procedure as mentioned above during the ramp-up 
phase~\cite{SHINANO2013};
See also Section \ref{sect:completeEnumeration}.
Until now, we have succeeded to solve tai30a, tai35b, tai40b and sko42 from 
QAPLIB \cite{QAPLIB} as reported in Section~\ref{sect:ISM}. 
But the C++ code needs further improvements for larger scale QAP instances. 
 
\subsection{Small scale instances}

\label{sect:MATLAB}

The main features of the parallel \BB method implemented in the MATLAB code 
can be summarized as follows:\vspace{-2mm} 
\begin{itemize}
\item 
The enumeration tree described in Section~\ref{sect:enumerationTree}.\vspace{-2mm} 
\item 
The Lagrangian DNN relaxation of sub QAPs and a modified version of 
NewtBracket \cite{KIM2020}.
See Section~\ref{sect:lowerBound}.\vspace{-2mm} 
\item 
A simple rounding of approximate optimal solutions of DNN$^p_{\lambda}(F,L)$ 
$((F,L) \in \Pi_r\times \Pi_r, \ r \in \{0,1,\ldots,n\})$ for the upper bounding procedure.\vspace{-2mm} 
\item 
The three types of branching rules, M, P and D, 
which are described in Sections~\ref{sect:QAPmean}, \ref{sect:DNNprimal} and~\ref{sect:DNNdual}, 
respectively.\vspace{-2mm}
\item 
A parallel breadth first search with processing each node QAP$(F,L)$ to compute $\zeta^l(F,L)$ and $\zeta^u(F,L)$ by one core.\vspace{-2mm}
\end{itemize}

The optimal values of all small QAP 
instances solved are known. 
The initial 
incumbent objective value is set to be the optimal value + 1, and it is updated to the optimal value 
which is found by the upper bound procedure at some iteration. With this setting, 
a high quality incumbent objective value close to the optimal value is known at the beginning of 
the \BB method. As a result, the role of the upper bounding procedure becomes less important, and the 
breadth first search is reasonable for parallel computation. 

All the computations for numerical results reported in Table~1, 2 and~3  
were performed using MATLAB 2020b on iMac Pro with Intel Xeon W CPU (3.2 GHZ), 8 cores 
and 128 GB memory. 

Numerical results on nug20 $\sim$ nug28, tai20a $\sim$ tai30b and bur20a $\sim$ bur20c 
are shown in Tables~1, 2 and 3, respectively. We observe that: \vspace{-2mm}
\begin{itemize}
\item The number of nodes (=the number of sub QAPs processed) greatly depends on the branching rules M, P and D.
Overall, the rule D performs  better than the others. \vspace{-2mm}
\item In particular, for bur26b and bur26c instances, the number of nodes generated 
with  the branching rule D is much smaller than those 
generated by the others.  \vspace{-2mm}
\item But, for the largest size QAP, tai30b, the branching rule M generates the smallest number of nodes. \vspace{-2mm} 
\end{itemize} 
We describe in detail 
on which branching rule is preferable for a given larger QAP instance in Section~\ref{sect:sampling}. 


\begin{table}[htp]
\begin{center}
\caption{
Numerical results on small scale QAP instances --- 1. Br : Branching rules presented in 
Section~\ref{sect:branchingRule}. 
} 
\label{table:numResMATLAB1}
\vspace{2mm}
\begin{tabular}{|l|r|r|r|r|r|r|r|c|}
\hline
QAP &            &    & No. of nodes             &  Total execution & \multicolumn{3}{|c|}{Time(sec) for computing} & No. of CPU \\
\cline{6-8}
instance & Opt.val & Br & \multicolumn{1}{|l|}{generated} & time(sec) in para.&    LB    &    UB    & Br & cores used \\
\hline
    &        &  M &    757       &    6.78e2    & 4.21e3 & 3.13e1 & 3.38e0 & \\
nug20     &       2,570 &  P &    301       &    3.18e2    & 1.80e3 & 1.45e1 & 1.18e1 & 8\\
     &        &  D &    412       &    3.84e2    & 2.15e3 & 1.91e1 & 3.11e1 & \\
\hline
     &        &  M &    312       &    4.93e2    & 2.67e3 & 1.80e1 & 1.48e0 & \\
nug21     &       2,438 &  P &    200       &    3.46e2    & 1.73e3 & 9.72e0 & 1.04e1 & 8 \\
     &        &  D &    272       &    3.73e2    & 1.75e3 & 1.33e1 & 2.59e1 & \\
\hline
     &        &  M &    429       &    7.16e2    & 4.43e3 & 2.55e1 & 3.37e0 & \\
nug22     &       3,596 &  P &    250       &    4.33e2    & 1.88e3 & 1.27e1 & 1.80e1 & 8 \\
     &        &  D &    149       &    4.10e2    & 2.08e3 & 9.96e0 & 2.36e1 & \\
\hline
     &        &  M &   1,413       &    3.34e3    & 2.24e4 & 9.30e1 & 2.16e1 & \\
nug24     &       3,488 &  P &    742       &    1.99e3    & 1.28e4 & 5.02e1 & 1.06e2 & 8 \\
     &       &  D &    913       &    2.57e3    & 1.46e4 & 7.39e1 & 2.46e2 & \\
\hline
     &        &  M &   8,446       &    2.00e4    & 1.49e5 & 5.60e2 & 1.53e2 & \\
nug25     &       3,744 &  P &   3,004       &    8.37e3    & 5.75e4 & 2.14e2 & 5.57e2 & 8 \\
     &        &  D &   3,805       &    8.97e3    & 6.46e4 & 3.41e2 & 1.29e3 & \\
\hline
     &        &  M &   2,810       &    1.48e4    & 1.06e5 & 4.33e2 & 1.27e2 & \\
nug27     &       5,234 &  P &    979       &    5.26e3    & 3.29e4 & 1.32e2 & 2.78e2 & 8 \\
     &        &  D &   2,170       &    9.83e3    & 6.91e4 & 3.90e2 & 1.41e3 & \\
\hline
     &        &  M &  10,977       &    5.58e4    & 4.31e5 & 2.27e3 & 5.79e2 & \\
nug28     &       5,166 &  P &   4,236       &    2.25e4    & 1.68e5 & 8.40e2 & 1.57e3 & 8 \\
     &        &  D &   9,977       &    4.55e4    & 3.37e5 & 2.15e3 & 7.83e3 & \\
\hline
\end{tabular}
\end{center}
\end{table}

\begin{table}[htp]
\begin{center}
\caption{
Numerical results on small scale QAP instances --- 2. Br : Branching rules presented in 
Section~\ref{sect:branchingRule}. 
} 
\label{table:numResMATLAB2}
\vspace{2mm}
\begin{tabular}{|l|r|r|r|r|r|r|r|c|}
\hline
QAP &            &    & No. of nodes             & Total execution & \multicolumn{3}{|c|}{Time(sec) for computing} & No. of CPU \\
\cline{6-8}
instance & Opt.val & Br & \multicolumn{1}{|l|}{generated} & time(sec) in para.&    LB    &    UB    & Br & cores used \\
\hline
     &      &  M &   2,444       &    1.49e3    & 1.04e4 & 1.06e2 & 7.72e-1 & \\
tai20a     &     703,482 &  P &   2,107       &    1.42e3    & 9.18e3 & 9.23e1 & 7.30e1 & 8 \\
     &      &  D &   1,600       &    1.17e3    & 8.30e3 & 7.31e1 & 1.02e2 & \\
\hline
     &   &  M &    183       &    2.03e2    & 2.34e2 & 1.13e0 & 3.12e-1 & \\
tai20b     &  122,455,319 &  P &    183       &    2.01e2    & 2.32e2 & 1.19e0 & 2.22e0 & 8 \\
     &   &  D &    146       &    2.40e2    & 2.38e2 & 1.23e0 & 6.11e0 & \\
\hline
     &   &  M &    367       &    9.68e2    & 2.76e3 & 1.18e1 & 2.54e-1 & \\
tai25b     &  344,355,646 &  P &    367       &    1.08e3    & 3.20e3 & 1.31e1 & 2.41e1 & 8 \\
     &   &  D &    367       &    1.08e3    & 2.90e3 & 1.26e1 & 6.21e1 & \\
\hline
     &   &  M &    989       &    1.11e4    & 7.03e4 & 3.87e2 & 6.83e1 & \\
tai30b     &  637,117,113 &  P &   1,654       &    1.87e4    & 1.06e5 & 3.97e2 & 6.26e2 & 8 \\
     &   &  D &   1,150       &    1.80e4    & 9.74e4 & 4.01e2 & 1.51e3 & \\
\hline
\end{tabular}
\end{center}
\end{table}

\begin{table}[htp]
\begin{center}
\caption{
Numerical results on small scale QAP instances --- 3. Br : Branching rules presented in 
Section~\ref{sect:branchingRule}. 
} 
\label{table:numResMATLAB3}
\vspace{2mm}
\begin{tabular}{|l|r|r|r|r|r|r|r|c|}
\hline
QAP &            &    & No. of nodes             & Total execution & \multicolumn{3}{|c|}{Time(sec) for computing} & No. of CPU \\
\cline{6-8}
instance & Opt.val & Br & \multicolumn{1}{|l|}{generated} & time(sec) in para.&    LB    &    UB    & Br & cores used \\
\hline
     &     &  M &   2,976       &    4.21e3    & 2.03e4 & 2.46e1 & 2.13e1 & \\
bur26a     &    5,426,670 &  P &  11,401       &    1.70e4    & 9.24e4 & 7.69e1 & 4.00e2 & 8 \\
     &     &  D &   2,358       &    3.35e3    & 1.52e4 & 1.80e1 & 1.81e2 & \\
\hline
           &              &  M  & 144,228       &    4.62e4 & 2.46e5 & 1.02e3 & 1.32e2 & \\
\cline{3-8}
bur26b     &    3,817,852 &  P  & More than     & More than  &        &        &        & 8  \\
           &              &     & 544,719       & 3 days     &       &        &        & \\
\cline{3-8}
           &               &  D &   8,625       &    5.15e3    & 2.51e4 & 7.49e1 & 2.42e2 & \\
\hline
     &     &  M &  16,990       &    6.66e3    & 3.18e4 & 1.20e2 & 2.12e1 & \\
bur26c    & 5,426,795 &  P    &  66,428       &    8.84e4    & 5.08e5 & 3.82e2 & 1.78e3 & 8 \\
     &     &  D &   2,416       &    3.06e3    & 9.90e3 & 1.71e1 & 9.59e1 & \\
\hline
\end{tabular}
\end{center}
\end{table}

\subsection{Challenging large scale instances}

\label{sect:ISM}
We solved challenging large scale instances, nug30, tai30a, tai35b, tai40b and sko42 on the ISM 
(Institute of Statistical Mathematics) supercomputer HPE SGI 8600,
which is a liquid cooled, tray-based, high-density clustered computer system. The ISM supercomputer has 384 computing nodes and 
each node has two Intel Xeon Gold 6154 3.0GHz CPUs (36 cores) and 384GB of memory. 
The \textsc{LoadCoordinator}\ 
process is assigned to a CPU core and each \textsc{Solver}\ 
process is assigned to a CPU core. 
Therefore, the number of \textsc{Solver}{s} 
is less than  the number of cores by one. 
All of the instances were solved as a single job, though check-pointing 
procedures were performed every 30 minutes
(see \cite{SHINANO2014} about check-pointing and restart mechanism).
Table \ref{table:numResLargeScale} shows the computational results.
Note that tai30a and sko42 were solved to the optimality for the first time.

\begin{table}[htp]
	\begin{center}
		\caption{
			Numerical results on challenging large scale QAP instances. Branching rule used was M. 
		} 
		\label{table:numResLargeScale}
		\vspace{2mm}
		\begin{tabular}{|l|r|r|r|r|}
			\hline
			QAP        &               &   No. of nodes                                         & \multicolumn{1}{c|}{Total execution} & \multicolumn{1}{c|}{No. of CPU } \\
			instance & Opt.val &  \multicolumn{1}{|l|}{generated} & time(sec) in para.&  \multicolumn{1}{c|}{cores used}    \\
			\hline
			nug30     &       6,124 &     26,181       &    3.14e3    & 1,728   \\
			\hline
			tai30a     &       1,818,146  &     34,000,579       &    5.81e5    & 1,728  \\
			\hline
			tai35b     &       283,315,445 &     2,620,547       &    2.49e5    &1,728 \\
			\hline
			tai40b     &       637,250,948 &     278,465       &    1.05e5    & 1,728   \\
			\hline
			sko42     &       15,812 &   6,019,419       &    5.12e5    & 5,184   \\
			\hline
		\end{tabular}
	\end{center}
\end{table}



\section{Some additional techniques for further developments}

\label{sect:additionalTechniques}

\subsection{A complete enumeration tree  generated by the simple branching rule M}

\label{sect:completeEnumeration} 

A distinctive feature of the branching rule M presented in Section~\ref{sect:QAPmean} 
is the independence from any lower 
bounding and upper bounding procedures. 
Therefore, we could create a complete enumeration tree 
using this rule 
in advance to start the \BB method. We usually  
start with the single root node associated with the original problem to be solved 
even if it is implemented in parallel.  

At the beginning of the \BB method with the normal ramp-up, only 
one \textsc{Solver}\ 
is used, and all others are idle till it finishes processing 
the original problem at the root node. 
As the \BB method proceeds, the number of \textsc{Solver}{s} 
to solve 
subproblems increases gradually. In the early stage of the parallel execution of such a 
\BB method for large scale problems on a large number of \textsc{Solver}{s}, 
many of the \textsc{Solver}{s} 
are idle, and it would take a while till all \textsc{Solver}{s} 
start 
running. This is clearly inefficient.

If we use branching rule M for a 
large scale QAP instance, we can start the \BB method at any level of the enumeration 
tree. For example, consider a QAP instance with $n=50$. Then we can 
create $50 \times 49 \times 48 = 117,600$ sub QAPs with dimension 47 at the fourth 
level in the enumeration tree using branching rules M. We can start the \BB method 
from these sub QAPs. We will investigate on the performance of this idea. 
\\

\subsection{Simple estimation of the number of nodes in the enumeration tree} 

\label{sect:sampling}

We have proposed $3$ types of branching rules M, P and D in 
Section~\ref{sect:branchingRule}, and observed that 
the number of nodes generated depends not only on 
instances 
but also the branching rule employed. 
A branching rule sometimes results in much fewer nodes than other rules as observed in  bur26a, 
bur26b and bur26c instances in 
Table 3. When  the \BB method is applied to larger scale QAP 
instances, selecting a good branching rule becomes a more critical issue.
For the selection,
we describe a simple sampling method to 
estimate the number of nodes at each depth (level) of the enumeration tree under 
a certain uniformity assumption which 
is necessary for the simplicity of the method 
and also sufficient for a rough estimation.  

Let G$(\NC,\EC)$ denote the enumeration tree 
to be constructed by applying
the \BB method to  QAP~\eqref{eq:QAP}. 
Recall that each node of $\NC$ 
is identified as a sub QAP, QAP$(F,L)$ for some 
$(F,L) \in \Pi_r\times\Pi_r \ (r = 0,1,\ldots,n)$. 
The node set $\NC$ is subdivided 
according to the 
depth of the enumeration tree G$(\NC,\EC)$, depending on their location: 
$\NC_r = \NC \cap \big(\Pi_r\times\Pi_r\big)$ $(r=0,\ldots,n)$. 
Let $m_r = \ \| \NC_r \|$ $(r=0,1,\ldots,n)$, 
which is estimated 
by the sampling method described below. 

Let $r=0$. Obviously $m_0 = 1$ since QAP$(\emptyset,\emptyset)$ is the unique element of 
$\NC_0$, which corresponds to the root node of the enumeration tree. 
Hence $\hat{m}_0 = m_0 = 1$. With the application of the lower 
bounding and upper bounding procedures to QAP$(\emptyset,\emptyset)$, 
it will become clear whether the \BB method terminates or 
proceeds to the branching procedure for generating $n-r = n$ child nodes of QAP$(\emptyset,\emptyset)$.  
In the former case, $\hat{m}_q = m_q= 0$ $(q=1.\ldots,n)$ and we are done. For the latter case,  
 we have $m_1 = n$ active nodes which forms $\NC_1$  and $\hat{m}_1 = m_1 = n$. 
If the lower bounding and upper bounding procedures  are applied
to all of the $\hat{m}_1 = n$ nodes, then the exact $m_2$ can be obtained; hence $\hat{m}_2 = m_2$. 
Instead, we choose $0 < s_1 < \hat{m}_1$ nodes randomly from the $\hat{m}_1$ nodes,  
and apply the procedures to those $s_1$ nodes to estimate $m_2$. As a result, 
some of the $s_1$ nodes are terminated. Let $t_1$ denote the number of the  nodes that remained active 
among the $s_1$ nodes. Applying the branching procedure to those $t_1$ nodes  
provides $(n-1)t_1$ nodes (included in $\NC_2$) 
as their child nodes. Thus, each node of the $s_1$ nodes  chosen randomly from the 
$\hat{m}_1 = m_1$ nodes have generated $(n-1)t_1/s_1$ child nodes on average. Consequently, the number $m_2$ 
of the child nodes from the $\hat{m}_1$ nodes is estimated by $\hat{m}_2 = \hat{m}_1 \times (n-1)t_1/s_1$. 

Assume that $\hat{m}_q$ $(q=0,1,\ldots,r)$ have been computed.
We describe how to compute
 $\hat{m}_{r+1}$  at the $r$th iteration of the sampling method. 
If $t_{r-1}$ attained $0$, {\it i.e.},  there are no active nodes 
from previous iteration, 
$\hat{m}_r$ was set to zero, so that we set $\hat{m}_{r+1} = 0$.  
Now we assume on the contrary that $0 < t_{r-1}$ nodes remained active,  
and $(n-(r-1))t_{r-1}$ nodes included in $\NC_r$ were generated  in the previous iteration. 
To compute an estimation $\hat{m}_{r+1}$ of $m_{r+1}$, choose $s_r$ nodes randomly 
from the $(n-(r-1))t_{r-1}$ nodes, and apply the 
lower bounding and upper bounding procedures to them. 
Suppose that there exist  $t_r$ active nodes. Then the branching procedure applied to them 
yields  $(n-r)t_r$ nodes in $\NC_{r+1}$. Therefore we set $\hat{m}_{r+1} = \hat{m}_{r} (n-r)t_r/s_r$. 

\begin{table}[htp]
\scriptsize{
\begin{center} 
\caption{
Estimation on the number of nodes in the enumeration tree for 
QAP instances from QAPLIB. 
'Estimated' $=$ $\sum_{r=0}^n \hat{m}_r$. 'Generated' $=$ $\sum_{r=0}^n m_r$.
All $(n-r)!$ feasible solutions of QAP$(F,L)$ were enumerated to compute its exact optimal value 
when $|F^c|=|L^c|=n-r = 7$ holds; hence $\hat{m}_{r} = m_r = 0$ for $r = n-6,n-5,\ldots,n$.
The \textcolor{blue}{\bf bold blue font} indicates the smallest estimate among the three branching rules M, P and D. 
- : not computed. (c) --- solved at ISM by the C++ code. (m) --- solved by the MATLAB code.
} 
\label{table:summary}
\vspace{2mm}
\begin{tabular}{|l||r|r||r|r||r|r|}
\hline 
               &  \multicolumn{6}{|c|}{Branching Rule (The number of nodes)}\\  
\hline 
               & \multicolumn{2}{|c||}{M } & \multicolumn{2}{|c||}{P } 
               & \multicolumn{2}{|c|}{D } \\ 
\hline                
Problem & Estimated  
& Generated & Estimated & Generated & Estimated & Generated \\ 
\hline
bur26a & 3,038 &  2,976     & 15,495 &  11,401      & \blue{\bf 2,404} &  2,358(m) \\ 
\hline
bur26b & 276,511 &  144,228   & 1,106,174 &  More than  & \blue{\bf 7,206} &  8,635(m) \\ 
        &        &            &           &  544,719    &              &           \\
\hline
bur26c & 18,114 &  16,990     & 102,053 &  66,428 & \blue{\bf{2,609}} &  2,416(m)  \\ 
bur26d & 119,212 &  -     & 7,583,919 &  -       & \blue{\bf 5,473} &  17,495(m)  \\ 
bur26e & 8,156 &  -     & 15,042 &  -       & \blue{\bf 1,621} &  1,621(m)  \\ 
bur26f & 432,056 &  -     & 6,530,633 &  -       & \blue{\bf 4,742} & 6,536(m) \\ 
bur26g & 32,935 &  -     & 18,023 &  -       & \blue{\bf 3,430} & 3,265(m)  \\ 
bur26h & 33,963 &  -     & 617,498 &  -       & \blue{\bf 3,645} & 4,211(m) \\ 
\hline
nug25  & 3239 &  8446(c)     &   \blue{\bf 2952} & 3004(m)  & 3,614 & 3,805(m) \\ 
nug27 &  2660 &  2610(c)     &  \blue{\bf 746}    &  979(m) & 2,272   & 2,170(m) \\ 
nug28 &  7,369 & 11,228(c) &   \blue{\bf 3,208} & 4,236(m) & 10,128 &  9,977(m) \\ 
nug30 &  12,419 & 26,297(c)&  10,905 & - & \blue{\bf 10,401} & - \\ 
\hline
tai30a & 42,472,865 & 34,000,579(c) &  19,280,014 & - & \blue{\bf 14,999,704}  & - \\ 
tai30b & 611               & 989(m)     & 611 & 1,654(m)  & 611 & 1,150(m)  \\ 
tai35a & 2.1e11 
& - &   3.9e10 
& - & \blue{\bf 1.3e10} 
& - \\ 
tai35b & 5,840,218 & 2,620,547(c)        &    4,618,822 & - & \blue{\bf 360,772} & - \\ 
tai40b &  212,954 & 278,465(c)              &   \blue{\bf 61,065} & - & 71,903,758    & - \\ 
tai50b & \blue{\bf 54,547,664} & -  &  1.4e9 
& - & 2.1e11 
& - \\ 
\hline
tho40 & 7.1e11 
&   -    & \blue{\bf 37,382,192} &-  & 97,834,698  & - \\ 
\hline
sko42 &  3,575,067 & 6,019,419(c)      &  3,071,294 &-  & \blue{\bf 2,764,606} & - \\ 
sko49 &  \blue{\bf{1.6e8}} 
& -   & 1.3e9 
&  - & 1.2e9 
& - \\ 
wil50  &  32,963,810               &  -   &   \blue{\bf 20,523,849} & - & 39,330,160 & - \\ 
\hline 
\end{tabular}
\end{center}
}
\end{table}

Table 5 shows the estimation of the total number of nodes for some instances. 
We observe that \vspace{-2mm}
\begin{description}
\item{(i) } (the total number of nodes generated in the numerical experiment)
\\ \mbox{ \ } \hspace{50mm} $\leq$ $3 \ \times $  (the estimate of the total number of nodes).\vspace{-2mm}
\item{(ii) } The branching rule P and D are expected to generate less nodes 
than the branching rule M except for tai50b and sko49.  
\vspace{-2mm}
\item{ (iii) } tai50b (using M), tho40 (P) and wil50 (P) could be solved.  
\vspace{-2mm}
\item{(iv) } tai35a is too difficult to solve using any of M, P and D.  \vspace{-2mm}
\item{(v) } We may challenge sko49 (M). 
\end{description}



\bibliographystyle{plain}
\bibliography{enhFOM}

\end{document}